\documentclass[12pt]{amsart}
\usepackage[a4paper,twoside,left=2cm,right=2cm,top=3.1cm,bottom=2.3cm]{geometry}

\usepackage{amsmath,amssymb,amscd,amsthm}
\usepackage{latexsym}
\usepackage{graphicx}
\usepackage[english]{babel}
\usepackage[latin1]{inputenc}       
\usepackage{textcomp}
\usepackage{times}
\usepackage{colortbl}

\newtheorem{theorem}{Theorem}[section]

\def\irr#1{{\rm Irr}(#1)}
\def\irrr#1#2 {\irr {#1 \mid #2}}

\newcommand{\R}{\mathbb R}

\begin{document}

\title[]{The strong log-concavity for first eigenfunction of the Ornstein-Uhlenbeck operator in the class of convex bodies}

\author[L. Qin]{Lei Qin }
\address{Institute of Mathematics,
Hunan University, Changsha, 410082, China}
\email{qlhnumath@hnu.edu.cn}


\maketitle

\baselineskip18pt

\parskip3pt

\begin{abstract}
In this paper, we prove that the first (positive) Dirichlet eigenvalue of the Ornstein-Uhlenbeck operator
\[
L(u)=\Delta u-(\nabla u,x),
\]
is strongly log-concave if the domain is bounded and convex, which improves the conclusion in \cite{Andrea-Paolo2024}. We also provide a characterization of the equality case of the Brunn-Minkowski inequality for the principal frequency of $L(u)$ in the class of convex bodies.
\end{abstract}

\section{Introduction}

A class of profoundly significant problems involves the study of the convexity properties of solutions to partial differential equations in the class of convex domains. There are many new techniques appeared and developed along these years (microscopic technique and macroscopic technique). A classical reference text is the book by Kawohl \cite{Kawohl1985}. Particularly, related to this paper are \cite{Korevaar-Lewis}, and \cite{Andrea-Paolo2024,Andrea-Lei-Paolo2024,Lei2}, where you can find more relevant references. Meanwhile, the Brunn-Minkowski type inequalities for variational functionals have a deep impact on both geometry and analysis. Examples are: Torsional rigidity, the first eigenvalues of Dirichlet Laplacian, Monge-Amp\`ere operator and so on, see for instance \cite{Acker1981,BL76,Caffarelli-Friedman,Colesanti2005,Korevaar-Lewis,Liu-Ma-Xu2010,Ma-Xu2008,Salani2005,Salani2012} and reference therein.

Very recently, a Brunn-Minkowski type inequality for the principal frequency of the Laplace operator in the Gauss space and the log-concavity of corresponding solution have been established, see  \cite{Andrea-Paolo2024}. Roughly speaking, the main goal of this paper is to establish the strong log-concavity of the solution to Gaussian principal frequency in the class of convex bodies, which relies on {\it Constant Rank theorem}, that is the rank of the Hessian of the corresponding convex solution is constant.  Then we establish the equality case of corresponding Brunn-Minkowski inequality in the class of convex bodies. In particular, based on a beautiful combination of the constant rank theorem and a continuity argument, we also deduce the main conclusion in \cite{Andrea-Paolo2024}, where authors prove the log-concavity of the Gaussian principal frequency by using an {\it infimal convolution} technique.

Here, we are concerned with the strong log-concavity of the first eigenfunction of the Ornstein-Uhlenbeck operator in the class of convex bodies. The log-concavity has been solved in \cite{Andrea-Paolo2024} in the class of convex bodies. The strong log-concavity has been established with extra assumptions, that is, the convex body is of class $C^{2,\alpha}_{+}$, for some $\alpha \in (0,1)$ and the domain is symmetric with respect to the origin. We note that as shown in \cite{Andrea-Paolo2024}, if the origin belongs to the domain and it is the maximum point of $u$, it is obvious that $(x,\nabla w(x))\geq 0$ where $w=-\ln u$ in $\Omega$, and then, we have $\Delta w>0$ in $\Omega$, following those steps in \cite{Caffarelli-Friedman,Korevaar-Lewis}, the strong log-concavity of function $u$ is obtained. In this paper, we remove the extra assumptions in \cite{Andrea-Paolo2024} and fill this gap.

Let $\Omega$ be an open, bounded Lipschitz domain. The Gaussian version of the principal frequency for $\Omega$ is defined by the following minimization problems:
\[
\lambda_{\gamma}:=\inf\left\{
\frac{\int_{\Omega} |\nabla u|^{2} d\gamma(x)}{\int_\Omega |u|^{2} d\gamma(x)}\colon u\in W^{1,2}_0(\Omega,\gamma),\ u\not\equiv0
\right\}.
\]
Here $\gamma$ denotes the Gaussian probability measure in $\R^n$, given by
$$
\gamma(A)= \frac{1}{(2\pi)^{n/2}} \int_{A} e^{-|x|^2/2} dx\,,\quad\text{for any measurable set }A\,.
$$
A standard variational method shows that there exists a (nontrivial) positive solution such that
\[
\lambda_{\gamma}(\Omega)=\int_\Omega|\nabla u|^2 d\gamma.
\]
And the solution is unique, up to multiplication by a positive constant. The standard regularity theory in \cite{Evans,GT} shows that $u\in C^{\infty}(\Omega)\cap C^0(\overline{\Omega})$ and $u$ satisfies the following boundary value problem
\begin{equation}\label{PDE1}
\left\{
\begin{aligned}
& -\Delta u +(x,\nabla u)=\lambda_\gamma u, \ \ \ \ \ \  \mathrm {in} \ \ \  \Omega,  \\
& u=0,  \ \ \ \ \ \  \ \ \ \ \ \ \ \ \  \mathrm {on} \ \ \ \partial \Omega.
\end{aligned}
\right.
\end{equation}
We refer to \cite{CCLaMP,Franceschi2024,H-Galyna2024} for a general study of the Gaussian principal frequency on various contexts. As for the $p$-Laplacian extension of the Ornstein-Uhlenbeck operator, please refer to \cite{Andrea-Lei-Paolo2024,Lei2}.

Our main results is as follows.

\begin{theorem}\label{T1}
Let $\Omega$ be an open, bounded and convex subset of $\R^n$, and let $u>0$ be a solution of problem \eqref{PDE1} in $\Omega$. Then the function
\[
w=-\ln u
\]
is strongly convex in $\Omega$, i.e.,
\[
D^2 w(x)>0, \quad \forall x\in \Omega.
\]
\end{theorem}

We note that equality case of the Brunn-Minkowski type inequality is typically associated with the strong log-concavity in the interior of the domain. Based on the above theorem, we can deduce the following theorem, which also improves the another result in \cite{Andrea-Paolo2024}. The main proof is very similar to the method in \cite{Colesanti2005,Andrea-Paolo2024}, we omit the proof.
\begin{theorem}
Let $\Omega_0,\Omega_1$ be an open, bounded and convex subsets of $\R^n$. Let $u_0$ and $u_1$ be solutions of \eqref{PDE1} with $\Omega=\Omega_0$ and $\Omega=\Omega_1$, respectively. Let $t\in [0,1]$ and set
\[
\Omega_t=(1-t)\Omega_0+t\Omega_1:=\{(1-t)x_0+tx_1:\ x_0\in \Omega_0,\ x_1\in \Omega_1 \}.
\]
If for some $t\in[0,1]$, equality holds, i.e.,
\[
\lambda_\gamma(\Omega_t)=(1-t)\lambda_\gamma(\Omega_0)+t\lambda_\gamma(\Omega_1),
\]
then $\Omega_0$ and $\Omega_1$ coincide up to translation.
\end{theorem}

\section{Preliminaries}

In this section, we list some preliminaries.

Let $\R^n$ be the $n$-dimensional Euclidean space. We say that the boundary of an open set $\Omega \subset \R^n$ is of class $C^{2,\alpha}_{+}$, for some $\alpha\in (0,1)$, if it is of class $C^{2,\alpha}$ and the Gauss curvature is strictly positive at every point of $\partial \Omega$. Let $u\in C^2(\Omega)$, we say that $u$ is strongly convex if the Hessian matrix $D^2 u$ is positive definite in $\Omega$; If the Hessian matrix $D^2 u$ is positive semi-definite in $\Omega$, we say that $u$ is convex in $\Omega$. Similarly, we say that a positive function $u\in C^2(\Omega)$ is strongly log-concave if the Hessian matrix $D^2 (\ln u)$ is negative definite in $\Omega$; If the Hessian matrix $D^2 (\ln u)$ is negative semi-definite in $\Omega$, we say that $u$ is log-concave in $\Omega$.

We denote the Gauss probability space by $( \mathbb{R}^{n}, \gamma)$, the measure $\gamma$ is given by
$$
\gamma(\Omega)= \frac{1}{(2\pi)^{n/2}} \int_{\Omega} e^{-|x|^2/2} dx,
$$
for any measurable set $\Omega\subseteq \mathbb{R}^{n}$. Throughout the paper, $d\gamma$ stands for integration with respect to $\gamma$, i.e., $d\gamma(x)=(2\pi)^{-n/2} e^{-|x|^2/2}dx$, and $d\gamma_{\partial \Omega}$ is the $(n-1)$-dimensional Hausdorff measure on $\partial\Omega$ with respect to $\gamma$. Let $\Omega\subseteq \R^n$ be a measurable set. We define $L^2(\Omega,\gamma)$ as follows:
\begin{equation}\nonumber
L^2(\Omega,\gamma)=\left\{ u:\Omega\rightarrow \R:\ \|u\|^{2}_{2,\gamma}:=\int_{\Omega} |u(x)|^2 d\gamma <+\infty \right\}.
\end{equation}
Similarly, $W^{1,2}_{0}(\Omega,\gamma)$ denotes the corresponding $\gamma$-weighted Sobolev space. Those Sobolev spaces are all separable Hilbert space with respect to Gaussian measure, (see for instance \cite{Tero1994,Radulescu2019} and references therein).

Let $\Omega$ be an open subset of $\R^n$, and $u\in C^2(\Omega)$. The Ornstein-Uhlenbeck operator is defined by
\[
Lu(x)=\Delta u(x)-(\nabla u(x),x).
\]
Integration by parts shows that $L$ is self-adjoint with respect to the Gaussian measure. If $\Omega$ is an open subsect of $\R^n$, and $\phi,\psi\in C_c^{\infty}(\Omega)$, then
\[
\int_\Omega \varphi L\psi d\gamma=-\int_\Omega (\nabla \varphi, \nabla \psi)d\gamma=\int_\Omega  \psi L\varphi d\gamma.
\]

In order to better understand the work of this paper, we list the main results in \cite{Andrea-Paolo2024}.
\begin{theorem}[Theorem 1.7 in \cite{Andrea-Paolo2024}]
Let $\Omega$ be an open bounded and convex subset of $\R^n$. Let $u$ be an eigenfunction of $\Omega$, that is a solution of problem \eqref{PDE1} in $\Omega$. Then the function is log-concave in $\Omega$.
\end{theorem}

\begin{theorem}[Theorem 1.2 in \cite{Andrea-Paolo2024}]
Let $\Omega_0$, $\Omega_1$ be open, bounded and Lipschitz subsets of $\R^n$. Let $t\in[0,1]$ and set
\[
\Omega_t=(1-t)\Omega_0+t\Omega_1.
\]
Then
\[
\lambda_\gamma(\Omega_t)\leq (1-t)\lambda_\gamma(\Omega_0)+t\lambda_\gamma(\Omega_1).
\]
\end{theorem}

\begin{theorem}[Theorem 1.3 in \cite{Andrea-Paolo2024}]
Let $\Omega_0$, $\Omega_1$ be open, bounded and convex subsets of $\R^n$. Assume moreover that $\Omega_0$ and $\Omega_1$ are symmetric with respect to the origin, and $\partial \Omega_0$ and $\partial \Omega_1$ are of class $C^{2,\alpha}_{+}$ for some $\alpha\in (0,1)$. If for some $t\in[0,1]$
\[
\lambda_\gamma((1-t)\Omega_0+t\Omega_1)= (1-t)\lambda_\gamma(\Omega_0)+t\lambda_\gamma(\Omega_1),
\]
then $\Omega_0=\Omega_1$.
\end{theorem}

\section{The proof of theorem \ref{T1}}
In this section, we give the proof of Theorem \ref{T1}.

Let $w=-\ln u$, we write the equation \eqref{PDE1} as follows:
\begin{equation}\label{PDE2}
\left\{
\begin{aligned}
& \Delta w =\lambda_\gamma+|\nabla w|^2+(x,\nabla w), \ \ \ \ \ \  \mathrm {in} \ \ \  \Omega,  \\
& w(x)\rightarrow +\infty,  \ \ \ \ \ \  \ \ \ \ \ \ \ \ \  \mathrm {as} \ \ x\rightarrow \ \partial \Omega.
\end{aligned}
\right.
\end{equation}
By the result in \cite{Andrea-Paolo2024}, we know that $w$ is convex in $\Omega$.

\begin{proof}[\bf Proof of Theorem \ref{T1}]

We divide its proof into two steps.

{\bf Claim 1:  } We prove that the matrix $D^2 w$ has a constant rank in $\Omega$.

Since $w$ is convex in $\Omega$, we know that all the eigenvalues of $D^2 w(x)$ are non-negative in $\Omega$.
Let $x_0$ be a point where the rank $r$ of $D^2 w$ is minimum. If $1\leq r\leq n$, the case has been considered in \cite[Theorem 4.5]{Andrea-Paolo2024}. We only need to consider the case $r=0$, that is
\[
\text{rank}( D^2 w(x_0))= 0.
\]
For $x\in \Omega$, we denote by $\varphi(x)$ the trace of the matrix $D^2 w$. Then
\[
\varphi(x)\geq 0\quad \text{in}\ \Omega\quad \text{and}\quad \varphi(x_0)=0.
\]
We will prove that there exist positive constants $c_1$ and $c_2$, such that
\begin{equation}\label{3-1}
\Delta \varphi(x)\leq c_1 |\nabla \varphi (x)|+c_2\varphi(x),\quad \text{in\ a neighborhood of}\ x_0.
\end{equation}
By the strong minimum principle, this implies that $\varphi\equiv 0$ in a neighborhood of $x_0$. Therefore, we conclude that the Hessian matrix $D^2 w$ has a constant rank in $\Omega$.

In the following, we prove inequality \eqref{3-1}. Let $\mathcal{U}\subseteq \Omega$ be a neighborhood of $x_0$. Let $z\in \mathcal{U}$, we may choose a coordinate system such that $D^2w(z)$ is a diagonal matrix (since the equation \eqref{PDE2} verified by $w$ is invariant under rotations). Then
\[
w_{ii}(z)\geq 0,\ 1\leq i\leq n,\quad w_{ij}(z)=0,\ 1\leq i\neq j\leq n.
\]
 By a direct calculation, we have
\[
\Delta \varphi=\mathop{\sum}_{i=1}^n (\Delta w)_{ii}\quad \sum_{i=1}^{n} (|\nabla w|^2 )_{ii}=2(\nabla w,\nabla \varphi)+2\sum_{i=1}^{n}|\nabla w_i|^2,
\]
and
\[
\sum_{i=1}^{n}(z,\nabla w)_{ii}=\sum_{i=1}^{n}2w_{ii}+(z,\nabla w_{ii})=2\varphi+(z,\nabla \varphi).
\]
Since $D^2 w$ is a diagonal matrix at point $z$, by equation \eqref{PDE2}, we have
\[
\Delta \varphi(z)=2\varphi+(z,\nabla \varphi)+2(\nabla w,\nabla \varphi)+2\sum_{i=1}^{n}w^2_{ii}.
\]
Since $\Omega$ is bounded, there exists a positive constant $c_1:=c_1(\Omega)$ such that $|x|\leq c_1$ in $\Omega$. Since $w$ is smooth in $\overline{\mathcal{U}}$, there exists a positive constant $c_2$, which depends on neighborhood $\mathcal{U}$ and $x_0$ such that
\[
|\nabla w|\leq c_2\quad \text{and}\quad |w_{ii}|\leq c_2,\quad \forall\ 1\leq i\leq n \quad \text{in\ }\ \overline{\mathcal{U}}.
\]
Thus, we have
\[
\Delta \varphi(z)  \leq (2c_1+c_2)|\nabla \varphi|+(2+c_2)\varphi.
\]
Since the constants can be chosen locally, we conclude that \eqref{3-1} hold.

{\bf Claim 2.} In the following, we claim that if the matrix $D^2 w$ has constant rank $r$ in $\Omega$, and $0\leq r\leq n-1$, we have $w$ is constant in $n-r$ coordinate directions or linear at least one direction.

If $r=0$, by the Claim 1, we have $D^2 w(x)$ is a zero matrix in $\Omega$. Differentialing equation \eqref{PDE2}, we deduce that $ w_i(x)=0$ in $\Omega$, for all $1\leq i\leq n$. Therefore, $w$ is constant in $n$ coordinate directions.

If $1\leq r\leq n-1$, following \cite[Proposition 4.5]{Andrea-Paolo2024} and \cite[Page 29-31]{Korevaar-Lewis}, we can conclude that through each point in $\Omega$, there is at least one line on which $w$ is linear.

Since $w\in C^{\infty}(\Omega)$, and $w(x)\rightarrow +\infty$, as $x\rightarrow x_0\in \partial \Omega$, by Claim 1 and Claim 2, we can conclude that the matrix $D^2 w$ has constant rank $n$ in $\Omega$, which implies that $w=-\ln u$ is strongly convex in $\Omega$. We complete the proof.
\end{proof}

\noindent
{\bf{Acknowledgements}} The author would like to thank Professors for their patient guidance and warm encouragement, and the Department of Mathematics and Computer Science of the University of Florence for the hospitality. The author also thanks China Scholarship Council (CSC) for the financial support during the visit at the University of Florence.

\bibliographystyle{amsplain}

\end{document}